\documentclass[12pt]{article}
\usepackage{}
\usepackage{mathrsfs}
\usepackage[justification=centering]{caption}
\usepackage{booktabs}
\usepackage{multirow}
\usepackage{array}
\usepackage[numbers,sort&compress]{natbib}
\usepackage[colorlinks,
            linkcolor=blue,
            anchorcolor=green,
            citecolor=magenta
            ]{hyperref}
\usepackage{mathrsfs,amssymb,amsfonts}
\usepackage{graphicx}
\usepackage{amsmath,amssymb,latexsym,color}
\usepackage[mathscr]{eucal}
\usepackage[numbers]{natbib}
\newtheorem{thm}{Theorem}[section]
\newtheorem{definition}[thm]{Definition}

\newtheorem{lemma}[thm]{Lemma}

\newtheorem{corollary}[thm]{Corollary}
\newtheorem{example}{Example}
\newtheorem{remark}{Remark}[section]

\newcommand{\proof}{{\it Proof.\quad}}
\newcommand{\qed}{\hfill\Box\medskip}
\usepackage{CJK}
\usepackage{setspace}

\begin{document}

\title{\bf The $h$-edge tolerable diagnosability of balanced hypercubes
\thanks{M. Xu's research is supported by the National Natural Science Foundation of China (11571044, 61373021) and the Fundamental Research Funds for the Central Universities.}}

\author{Min Xu\footnote{\footnotesize Corresponding author.
 {\em E-mail address:} xum@bnu.edu.cn (M. Xu).}\quad Yulong Wei\\
{\small  School of Mathematical Sciences, Beijing Normal University}\\
\noindent {\small Laboratory of Mathematics and Complex Systems,
Ministry of Education,}\\ \noindent {\small Beijing, 100875, China} }
 \date{}
\maketitle

\setlength{\baselineskip}{24pt}

\noindent {\bf Abstract}\quad To measure the fault diagnosis capability of a multiprocessor system with faulty links, Zhu et al. [Theoret. Comput. Sci. 758 (2019) 1--8] introduced the $h$-edge tolerable diagnosability. This kind of diagnosability is a generalization of the concept of traditional diagnosability. In this paper, as complement to the results in [Theoret. Comput. Sci. 760 (2019) 1--14], we completely determine the $h$-edge tolerable diagnosability of balanced hypercubes $BH_n$ under the PMC model and the MM$^*$ model. Thus, the traditional diagnosability of $BH_n$ is also determined.

\noindent {\bf Keywords}\quad Balanced hypercube; Fault diagnosis; PMC model; MM$^*$ model

\vskip0.6cm

\section{Introduction}
Processor failure has become an ineluctable event in a large-scale multiprocessor system. Thus, to keep the multiprocessor system performing its functions efficiently and economically, recognizing faulty processors correctly is a task of top priority. The process of recognizing faulty processors
in a multiprocessor system is called {\em fault diagnosis}, and the {\em diagnosability} of a system is the maximum number of faulty processors the system can recognize. Historically, many scholars and researchers proposed different models to investigate fault diagnosis. In 1967, Preparata, Metze and Chien \cite{Pre} proposed the PMC model for fault diagnosis in multiprocessor systems. Under this model, all adjacent processors of a system can test one another.
In 1992, by modifying the MM model \cite{Mae}, Sengupta and Dahbura \cite{Sen} proposed the MM$^*$ model assuming that each processor has to test two processors if the processor is adjacent to the latter two processors. Some references related to fault diagnosis studies under the PMC model or MM$^*$ model can be seen in \cite{CDH,FanL,Gu,Lai,Lin17,Pen,WanH,Wei,Wei18,Wei186,Wei188,XuTH,XuTZ,Yang,Yang13,YuanL,ZhaX,Zho,ZhuL,Zha}.

In recent literature \cite{ZhuL}, Zhu et al. introduced the {\em $h$-edge tolerable diagnosability} $t^e_h(G)$ to measure the fault diagnosis capability of a multiprocessor system $G$ with faulty links. This kind of diagnosability is a generalization of the concept of traditional diagnosability. Specifically, $t^e_h(G)$ is the minimum diagnosability of graphs $G-F_e$ which satisfy that $F_e\subseteq E(G)$ and $|F_e|\leq h$. Note that if a processor $u$ has no fault-free neighbors, it is impossible to determine whether $u$ is faulty or not in the fault diagnosis. Then $t^e_h(G)=0$ for $h\geq k$, where $G$ is a $k$-regular graph. Hence, a key issue for the $h$-edge tolerable diagnosability of a $k$-regular graph study is the case of $0\leq h\leq k$.

Let $C(G)$ be the maximum number of common neighbors of any two vertices in the graph $G$. Wei and Xu determined the $h$-edge tolerable diagnosabilities of regular graphs as follows.
\begin{thm}[\cite{Wei188}]\label{I1}
Let $G=(V, E)$ be a $k$-regular triangle-free graph with $k\geq2$. If $C(G)\leq k-1$, then $t_h^e(G)=k-h$ under the PMC model for $0\leq h \leq k$.
\end{thm}
\begin{thm}[\cite{Wei188}]\label{I2}
Let $G=(V, E)$ be a connected $k$-regular triangle-free graph with $k\geq3$. If $C(G)\leq k-1$, then
\begin{equation*}
  t_h^e(G)=
\left\{
  \begin{array}{lll}
   2 & \hbox{if $G$ is isomorphic to $G_8$ and $h=0$}; \\
   \\
   k-1 & \hbox{if $G$ is isomorphic to $G_{k+1, k+1}$ and $h=0$}; \\
    \\
    k-h & \hbox{otherwise}
  \end{array}
\right.
\end{equation*}
under the MM$^*$ model for $0\leq h \leq k$.
\end{thm}
In Theorem \ref{I2}, $G_8$ is the graph with vertex set $V(G_8)=\{x_1, x_2, \ldots, x_8\}$ and edge set $E(G_8)=\{x_ix_{i+1}\mid 1\leq i\leq7\}\cup\{x_8x_1\}\cup\{x_jx_{j+4}\mid 1\leq j\leq4\}$, and
$G_{k+1, k+1}$ is the graph with vertex set $V(G_{k+1, k+1})=\{x_1, x_2, \ldots, x_{k+1}, y_1, y_2, \ldots, y_{k+1}\}$ and edge set $E(G_{k+1, k+1})=\{x_iy_j\mid 1\leq i\leq k+1, 1\leq j\leq k+1, i\neq j\}$.

In this paper, we are concerned with the fault diagnosis capability analysis of balanced hypercubes $BH_n$. The $n$-dimensional balanced hypercube $BH_n$, as one of important variants of the well-known hypercubes, was proposed by Wu and Huang \cite{WuH}.
In recent years, $BH_n$ has received considerable attention. For example, Yang \cite{Yang,Yang13} studied the conditional diagnosability of $BH_n$ under the PMC model and the MM$^*$ model. Gu et al. \cite{Gu} determined the $1,2$-good-neighbor diagnosability of $BH_n$ under the PMC model and the MM$^*$ model. Lin et al. \cite{Lin17} determined the $1,2,3$-extra conditional fault-diagnosability of $BH_n$ under the PMC model. Zhang et al. \cite{ZhaX} investigated the $(t, k)$-diagnosability of $BH_n$ under the PMC model.
Although $BH_n$ is a $2n$-regular and triangle-free graph, $C(BH_n)=2n$. Thus, $BH_n$ does not satisfy the conditions of Theorems \ref{I1} and \ref{I2}. As complement to Theorems \ref{I1} and \ref{I2}, we establish the $h$-edge tolerable diagnosability of balanced hypercubes $BH_n$ under the PMC model and the MM$^*$ model for $0\leq h \leq 2n$ and $n\geq1$.
\medskip

\section{Terminology and preliminaries}\label{2}
A graph $G =\big(V(G), E(G)\big)$ is used to represent a system (or a network), where each vertex of $G$ represents a processor and each edge of $G$ represents a link. The {\em connectivity} $\kappa(G)$ is the minimum cardinality of all vertex subsets $S\subseteq V(G)$ satisfying that $G-S$ is disconnected or trivial. The {\em neighborhood} $N_G(v)$ of a vertex $v$ in $G$ is the set of vertices adjacent to $v$. We refer readers to \cite{Bon} for terminology and notation unless stated otherwise.

In 1997, Wu and Huang proposed balanced hypercubes $BH_n$. We restate the definition of $BH_n$ as follows.
\begin{definition}[\cite{WuH}]\label{D0}
The $n$-dimensional balanced hypercube $BH_n=(V(BH_n), E(BH_n))$ has vertex set $V(BH_n)=\{(a_0, a_1, \ldots, a_i, \ldots, a_{n-1})\mid a_i\in\{0, 1, 2, 3\}, 0\leq i\leq n-1\}$. Each vertex $(a_0, a_1, \ldots, a_{i-1}, a_i, a_{i+1}, \ldots, a_{n-1})$ of $BH_n$ has $2n$ neighbors:
\begin{description}
  \item[(1)] $((a_0\pm 1)~{\rm mod~4}, a_1, \ldots, a_{i-1}, a_i, a_{i+1},\ldots, a_{n-1})$,
  \item[(2)] $((a_0\pm 1)~{\rm mod~4}, a_1, \ldots, a_{i-1}, (a_i+(-1)^{a_0})~{\rm mod~4}, a_{i+1},\ldots, a_{n-1})$.
\end{description}
\end{definition}

\begin{figure}[hptb]
  \centering
  \includegraphics[width=8cm]{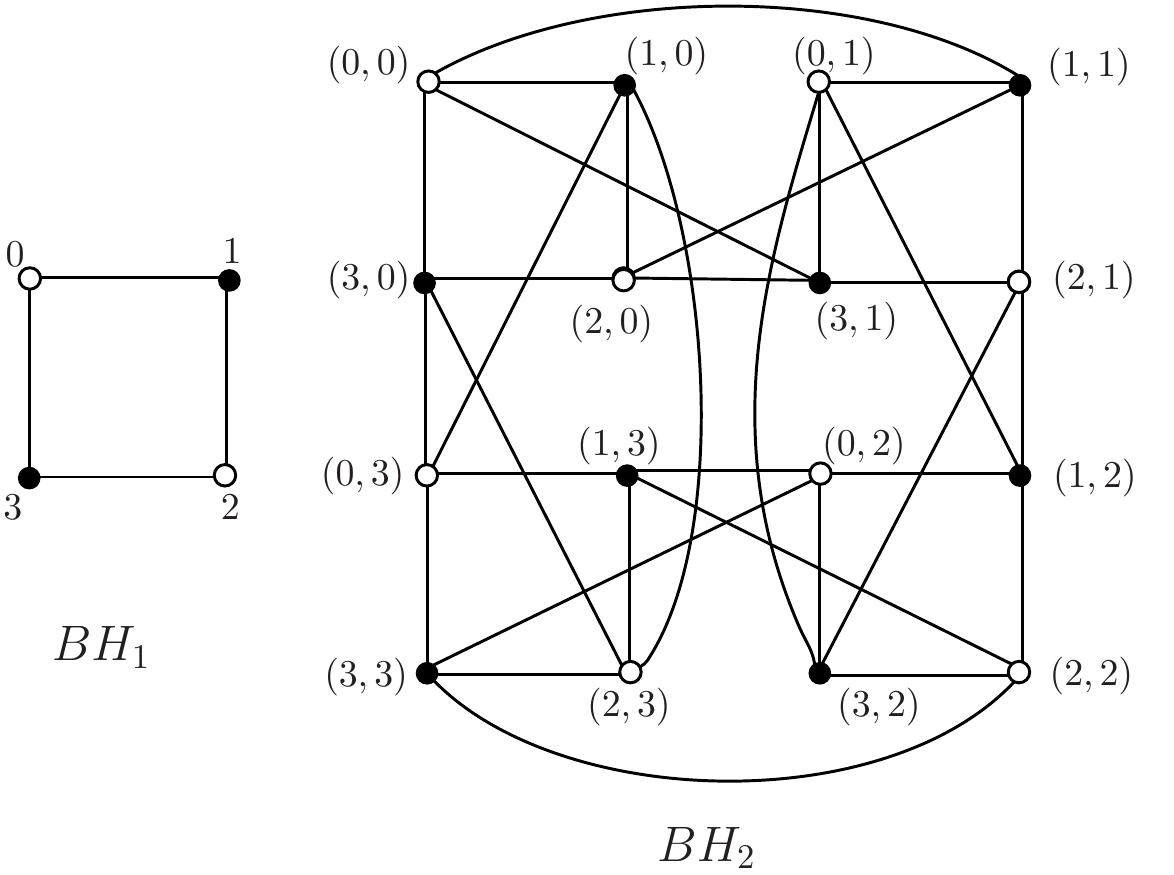}\\
  \caption{Illustration of $BH_1$ and $BH_2$. }\label{DD1}
\end{figure}

Figure \ref{DD1} shows $BH_1$ and $BH_2$.

Some basic but useful properties of $BH_n$ are presented as follows.

\begin{lemma}[\cite{WuH}]\label{pp4}
The balanced hypercube $BH_n$ is bipartite and $\kappa(BH_n)=2n$.
\end{lemma}

\begin{lemma}[\cite{Yang}]\label{pp8}
Let $u$ be an arbitrary vertex of $BH_n$ for $n\geq1$. Then, for an arbitrary vertex $v$ of $BH_n$, either $|N_{BH_n}(u)\cap N_{BH_n}(v)|=0$, $|N_{BH_n}(u)\cap N_{BH_n}(v)|=2$, or $|N_{BH_n}(u)\cap N_{BH_n}(v)|=2n$. Furthermore, there is exactly one vertex $w$ such that $|N_{BH_n}(u)\cap N_{BH_n}(w)|=2n$.
\end{lemma}

Now, we introduce the definition of the traditional diagnosability of a graph.

\begin{definition}[\cite{Dah}]\label{D3}
A graph $G=(V, E)$ of $n$ vertices is $t$-diagnosable if all faulty vertices can be detected without replacement, provided that the number of faults does not exceed $t$. The diagnosability $t(G)$ of a graph $G$ is the maximum value of $t$ such that $G$ is $t$-diagnosable.
\end{definition}

For any two sets $A$ and $B$, we use $A-B$ to denote a set obtained by removing all elements of $B$ from $A$. The {\em symmetric difference} of two sets $F_1$ and $F_2$ is defined as the set $F_1\bigtriangleup F_2$ $=(F_1-F_2)\cup (F_2-F_1)$. The following lemmas give necessary and sufficient conditions for a graph to be $t$-diagnosable under the PMC model and the MM$^*$ model.

\begin{figure}[hptb]
 \centering
  \includegraphics[width=10cm]{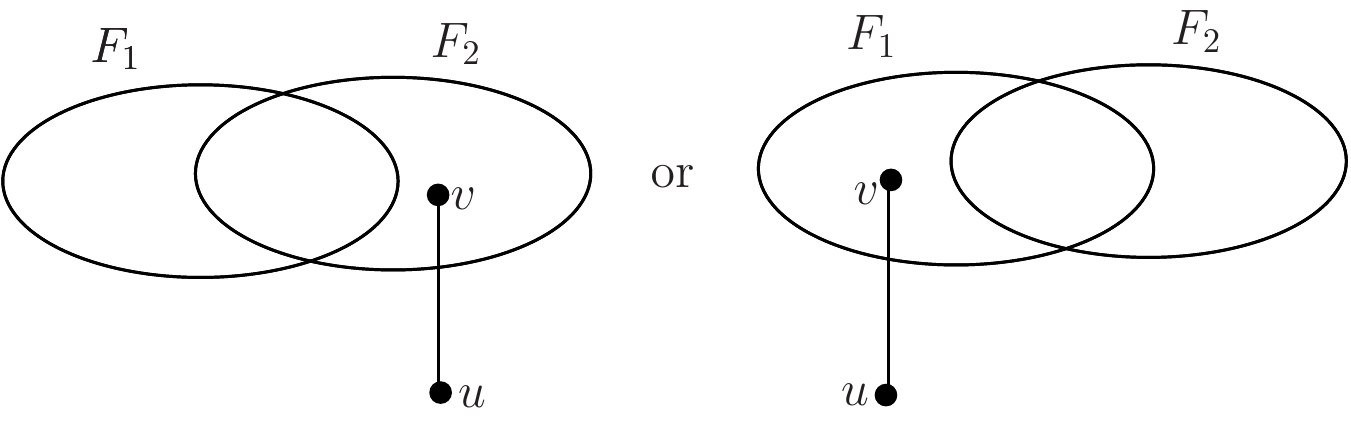}\\
  \caption{ The illustration of Lemma \ref{L01}.
}\label{L11}
\end{figure}

\begin{figure}[hptb]
  \centering
  \includegraphics[width=6cm]{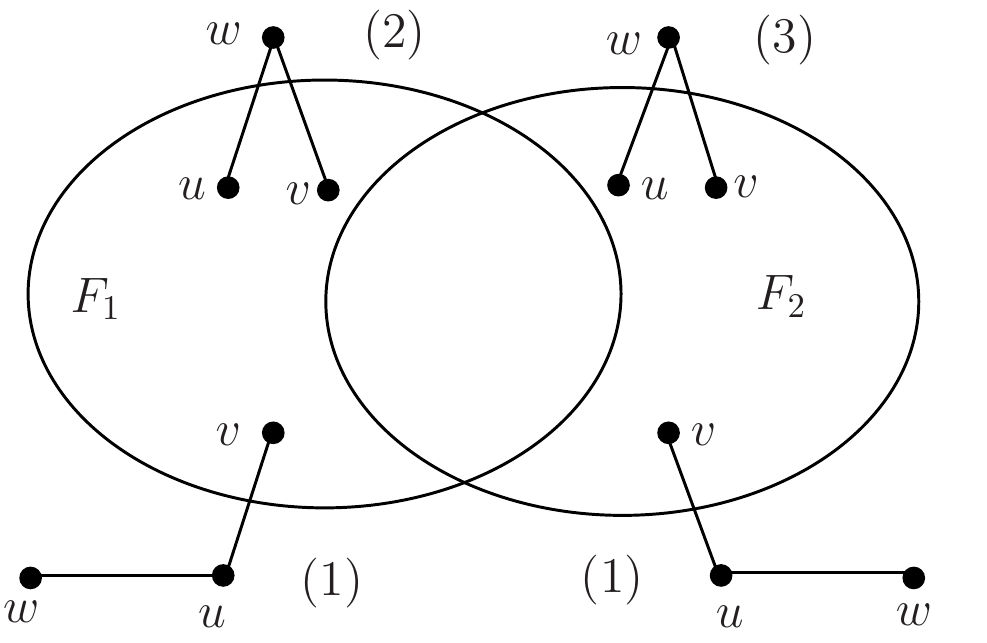}\\
  \caption{ The illustration of Lemma \ref{L02}.
}\label{L22}
\end{figure}

\begin{lemma}[\cite{Dah}]\label{L01}
A graph $G =(V, E)$ is $t$-diagnosable under the PMC model if and only if for any two distinct subsets $F_1$ and $F_2$ of $V$ with $|F_1|\leq t$ and $|F_2|\leq t$, there exists a test from $V-(F_1\cup F_2)$ to $F_1\bigtriangleup F_2$ {\rm(}see Figure \ref{L11} {\rm)}.
\end{lemma}
\begin{lemma}[\cite{Sen}]\label{L02}
A graph $G =(V, E)$ is $t$-diagnosable under the MM$^*$ model if and only if for any two distinct subsets $F_1$ and $F_2$ of $V$ with $|F_1|\leq t$ and $|F_2|\leq t$, at least one of the following conditions is satisfied {\rm(}see Figure \ref{L22} {\rm)}:
\begin{enumerate}
\item[{\rm (1)}] There are two vertices $u, w\in V-(F_1\cup F_2)$ and there is a vertex $v\in F_1\bigtriangleup F_2$
such that $uv\in E$ and $uw\in E$.

\item[{\rm (2)}] There are two vertices $u, v\in F_1-F_2$ and there is a vertex $w\in V-(F_1\cup F_2)$
such that $uw\in E$ and $vw\in E$.

\item[{\rm (3)}] There are two vertices $u, v\in F_2-F_1$ and there is a vertex $w\in V-(F_1\cup F_2)$
such that $uw\in E$ and $vw\in E$.
\end{enumerate}
\end{lemma}

We call sets $F_1$ and $F_2$ {\em distinguishable} under the PMC (resp. MM$^*$) model if they satisfy the condition of Lemma \ref{L01} (resp. at least one of the conditions of Lemma \ref{L02}); otherwise, $F_1$ and $F_2$ are said to be {\em indistinguishable}.

Recently, Zhu et al. introduced the definition of the $h$-edge tolerable diagnosability of graphs as follows.

\begin{definition}[\cite{ZhuL}]\label{D1}
Given a diagnosis model and a graph $G$, $G$ is $h$-edge tolerable $t$-diagnosable under the
diagnosis model if for any edge subset $F_e$ of $G$ with $|F_e|\leq h$, the graph $G-F_e$ is $t$-diagnosable under the
diagnosis model. The $h$-edge tolerable diagnosability of $G$, denoted as $t_h^e(G)$, is the maximum integer $t$ such that $G$ is $h$-edge tolerable $t$-diagnosable.
\end{definition}
Clearly, $t_0^e(G)=t(G)$ holds for any graph $G$.

\section{Main Results}\label{3}
In this section, we investigate the $h$-edge tolerable
diagnosability of a balanced hypercube $BH_n$ under the PMC model and the MM$^*$ model.

Wei and Xu gave an upper bound of the $h$-edge tolerable
diagnosability of a $k$-regular graph $G$ under the PMC model and the MM$^*$ model as follows.
\begin{lemma}[\cite{Wei188}]\label{LL1}
Let $G=(V, E)$ be a $k$-regular graph with $k\geq2$. Then $t_h^e(G)\leq k-h$ under the PMC model and the MM$^*$ model for $0\leq h \leq k$.
\end{lemma}

Thus, we immediately obtain the upper bound of the $h$-edge tolerable
diagnosability of a balanced hypercube $BH_n$ under the PMC model and the MM$^*$ model.
\begin{corollary}\label{LLU}
Let $BH_n$ be an $n$-dimensional balanced hypercube. Then $t_h^e(BH_n)\leq 2n-h$ under the PMC model and MM$^*$ model for $0\leq h \leq 2n$.
\end{corollary}

Now, we give a lower bound of the $h$-edge tolerable
diagnosability of a balanced hypercube $BH_n$ under the MM$^*$ model. In the following statements, for a vertex subset $A$ of a graph $G$, we use $N_G(A)$ to denote the set $\big(\bigcup _{v\in A}N_G(v)\big)-A$.

\begin{lemma}\label{LLL}
Let $BH_n$ be an $n$-dimensional balanced hypercube with $n\geq2$. Then $t_h^e(BH_n)\geq 2n-h$ under the MM$^*$ model for $0\leq h \leq 2n$.
\end{lemma}
\proof
For an arbitrary edge subset $F_e\subseteq E(BH_n)$ with $|F_e|\leq h$, suppose that there exist two distinct vertex subsets  $F_1, F_2\subseteq V(BH_n)$ such that $F_1$ and $F_2$ are indistinguishable in $BH_n-F_e$ under the MM$^*$ model. We will prove the lemma by showing that $|F_1|\geq 2n-h+1$ or $|F_2|\geq 2n-h+1$. If $|F_1\cap F_2|\geq 2n-h$, then $|F_1|\geq 2n-h+1$ or $|F_2|\geq 2n-h+1$. Now, we assume that $|F_1\cap F_2|\leq 2n-h-1$. Our discussion is divided into two cases as follows.

\medskip

{\em Case 1.} For each vertex $u\in F_1\bigtriangleup F_2$, $N_{BH_n-F_e}(u)-(F_1\cup F_2)=\emptyset$.

\begin{figure}[hptb]
  \centering
  \includegraphics[width=6cm]{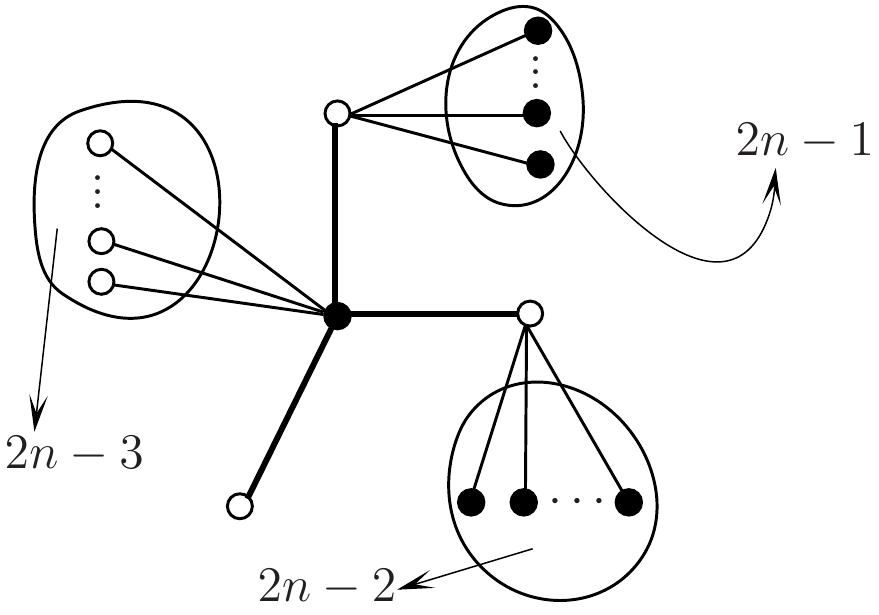}\\
  \caption{Illustration for a lower bound of $|N_{BH_n}(V(K_{1,3}))|$.}\label{C1}
\end{figure}

In this case, choose a vertex $x\in F_1\bigtriangleup F_2$. Then $N_{BH_n-F_e}(x)\subseteq F_1\cup F_2$ and $|N_{BH_n-F_e}(x)|\geq 2n-h>|F_1\cap F_2|$. Thus, there exists a vertex $y\in F_1\bigtriangleup F_2$ such that $xy\in E(BH_n-F_e)$. We have $N_{BH_n-F_e}(\{x, y\})\subseteq F_1\cup F_2$. By Lemma \ref{pp4}, we know that $BH_n$ is a bipartite graph. Note that $|N_{BH_n-F_e}(\{x, y\})|-|F_1\cap F_2|\geq(4n-2-h)-(2n-h-1)=2n-1\geq3$ for $n\geq2$. Then, there exists a star $K_{1,3}\subseteq BH_n[F_1\bigtriangleup F_2]-F_e$.
Note that $|N_{BH_n}(V(K_{1,3}))|\geq (2n-3)+(2n-1)+(2n-2)=6n-6$ (see Figure \ref{C1}).
Since $N_{BH_n-F_e}(V(K_{1,3}))\cup V(K_{1,3})\subseteq F_1\cup F_2$,
\begin{eqnarray*}
|F_1\cup F_2|&\geq&|N_{BH_n-F_e}(V(K_{1,3}))|+|V(K_{1,3})|\\
&\geq&(6n-6)-h+4\\
&=&(4n-2h+1)+(2n+h-3)\\
&\geq&4n-2h+1.
\end{eqnarray*}
The last inequality holds for $n\geq2$. Then $|F_1|\geq 2n-h+1$ or $|F_2|\geq 2n-h+1$.

\medskip

{\em Case 2.} $N_{BH_n-F_e}(F_1\bigtriangleup F_2)-(F_1\cup F_2)\neq\emptyset$.

Without loss of generality, suppose that $u\in F_1-F_2$, $v\in V(BH_n)-(F_1\cup F_2)$ such that $uv\in E(BH_n-F_e)$.
Note that $F_1$ and $F_2$ are indistinguishable in $BH_n-F_e$ under the MM$^*$ model.  Then $N_{BH_n-F_e}(v)-\{u\} \subseteq F_2$ and $|N_{BH_n-F_e}(v)\cap (F_2-F_1)|\leq1$. Thus, $|F_1\cap F_2|\geq |N_{BH_n-F_e}(v)|-2\geq 2n-h-2$. If $|F_1-F_2|\geq3$ or $|F_2-F_1|\geq3$, then $|F_1|\geq 2n-h+1$ or $|F_2|\geq 2n-h+1$.

Next, we suppose that $|F_1-F_2|\leq2$ and $|F_2-F_1|\leq2$.

{\em Case 2.1.} For each vertex $w\in F_1\bigtriangleup F_2$, $|N_{BH_n-F_e}(w)-(F_1\cup F_2)|\leq1$.

Note that $N_{BH_n-F_e}(\{u, v\})-(F_1-F_2)\subseteq F_2$. Thus,
\begin{eqnarray*}
|F_2|&\geq& |N_{BH_n-F_e}(\{u, v\})-(F_1-F_2)|\\
&\geq& (2n-2)+(2n-1)-h\\
&=& (2n-h+1)+(2n-4)\\
&\geq& 2n-h+1
\end{eqnarray*}
for $n\geq2$. The second inequality holds for $BH_n$ is a bipartite graph.

\medskip

{\em Case 2.2.} For some vertex $w\in F_1\bigtriangleup F_2$, $|N_{BH_n-F_e}(w)-(F_1\cup F_2)|\geq2$.

\begin{figure}[hptb]
  \centering
  \includegraphics[width=6cm]{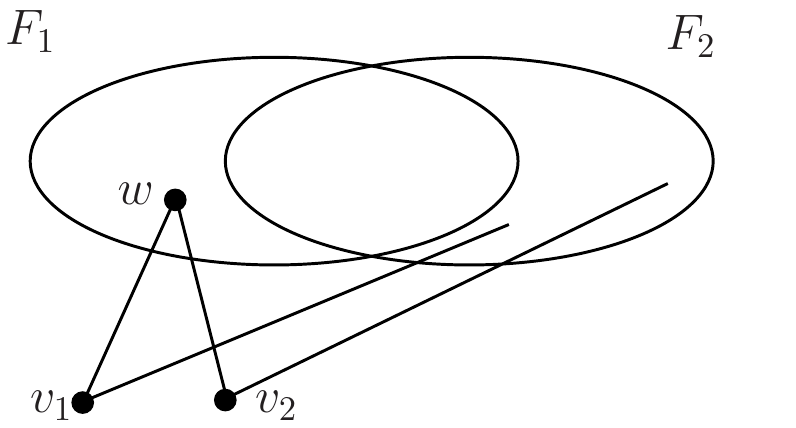}\\
  \caption{Illustration for $N_{BH_n}(v_1)\neq N_{BH_n}(v_2)$. }\label{f1}
\end{figure}

If there exists a vertex subset $\{v_1, v_2\}\subseteq N_{BH_n-F_e}(w)-(F_1\cup F_2)$ for some vertex $w\in F_1\bigtriangleup F_2$ such that $N_{BH_n}(v_1)\neq N_{BH_n}(v_2)$, then by Lemma \ref{pp8}, $|N_{BH_n}(v_1)\cap N_{BH_n}(v_2)|=2$. Without loss of generality, assume that $w\in F_1-F_2$. Since $F_1$ and $F_2$ are indistinguishable in $BH_n-F_e$ under the MM$^*$ model, $N_{BH_n-F_e}(v_1)-\{w\}\subseteq F_2$ and $N_{BH_n-F_e}(v_2)-\{w\}\subseteq F_2$ (see Figure \ref{f1}). Thus,
\begin{eqnarray*}
|F_2|&\geq& |N_{BH_n-F_e}(\{v_1, v_2\})-\{w\}|\\
&\geq& (2n-1)+(2n-1)-(2-1)-h\\
&=& (2n-h+1)+(2n-4)\\
&\geq& 2n-h+1.
\end{eqnarray*}
The last inequality holds for $n\geq2$.

\begin{figure}[hptb]
  \centering
  \includegraphics[width=6cm]{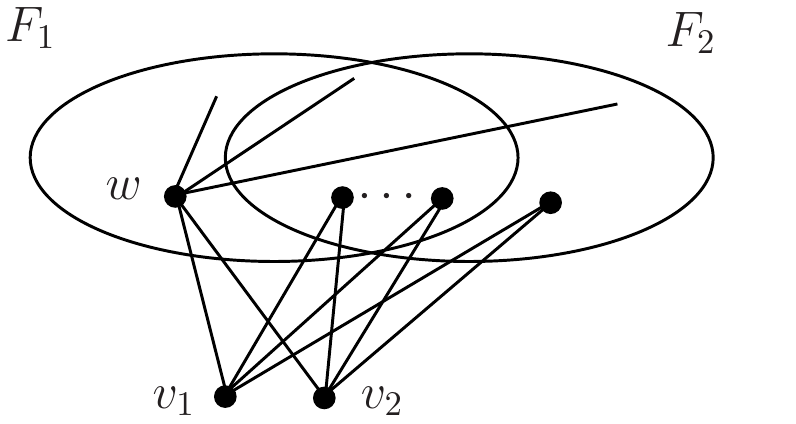}\\
  \caption{Illustration for $N_{BH_n}(v_1)=N_{BH_n}(v_2)$. }\label{f2}
\end{figure}

Otherwise, by Lemma \ref{pp8}, for each vertex $w\in F_1\bigtriangleup F_2$, $|N_{BH_n-F_e}(w)-(F_1\cup F_2)|\leq2$. Without loss of generality, assume that $w\in F_1-F_2$ and let $\{v_1, v_2\}=N_{BH_n-F_e}(w)-(F_1\cup F_2)$, where $N_{BH_n}(v_1)=N_{BH_n}(v_2)$ (see Figure \ref{f2}). Then $|F_1\cap F_2|\geq2n-2-\left\lfloor\dfrac{h}{2}\right\rfloor$ and $N_{BH_n-F_e}(w)-\{v_1, v_2\}\subseteq F_1\cup F_2$. Note that $(N_{BH_n-F_e}(w)-\{v_1, v_2\})\cap N_{BH_n-F_e}(v_1)=\emptyset$ and $|F_1\cup F_2|\geq |N_{BH_n-F_e}(\{w, v_1\})\cup\{w\}-\{v_2\}|\geq (2n-2)+(2n-1)+1-h=4n-h-2$.
Thus,
\begin{eqnarray*}
|F_1|+|F_2|&=& |F_1\cup F_2|+|F_1\cap F_2|\\
&\geq& (4n-h-2)+(2n-2-\left\lfloor\frac{h}{2}\right\rfloor)\\
&=&(4n-2h+1)+(\left\lceil\frac{h}{2}\right\rceil+2n-5).
\end{eqnarray*}
Note that $n\geq2$. If $n\geq3$ or $1\leq h\leq2n$, then $|F_1|+|F_2|\geq 4n-2h+1$ which means $|F_1|\geq 2n-h+1$ or $|F_2|\geq 2n-h+1$. Now, we assume that $n=2$ and $h=0$.

If $|F_1-F_2|=1$, then
\begin{eqnarray*}
|F_2|&\geq& |N_{BH_n-F_e}(\{w, v_1\})-\{v_2\}|\\
&\geq& (2n-2)+(2n-1)-h\\
&=& 2n-h+1
\end{eqnarray*}

If $|F_1-F_2|=2$ and $|F_1\cap F_2|\geq 2n-1-\left\lfloor\dfrac{h}{2}\right\rfloor$, then $|F_1|\geq2n-h+1$.

If $|F_1-F_2|=2$ and $|F_1\cap F_2|=2n-2-\left\lfloor\dfrac{h}{2}\right\rfloor=2$, then $BH_2-(F_1\cap F_2)$ is connected owing to $|F_1\cap F_2|=2<4=\kappa(BH_2)$ by Lemma \ref{pp4}. By the assumption that $|N_{BH_n-F_e}(w)-(F_1\cup F_2)|\leq2$ for each vertex $w$ in $F_1\bigtriangleup F_2$, we have $|N_{BH_2}(F_1\bigtriangleup F_2)-(F_1\cup F_2)|\leq 2\times|F_1\bigtriangleup F_2|\leq8$ and so $|V(BH_2)-(F_1\cup F_2)-N_{BH_2}(F_1\bigtriangleup F_2)|= |V(BH_2)|-(|F_1-F_2|+|F_1\cap F_2|+|F_2-F_1|)-|N_{BH_2}(F_1\bigtriangleup F_2)-(F_1\cup F_2)|\geq 2^4-(2+2+2)-8\geq2$.
Thus, there exists a vertex $x$ in $V(BH_2)-(F_1\cup F_2)$ with $x\notin N_{BH_2}(F_1\bigtriangleup F_2)$ connected to a vertex in $F_1\bigtriangleup F_2$ by a path of $BH_2-(F_1\cap F_2)$, which contradicts that $F_1$ and $F_2$ are indistinguishable in $BH_2$ under the MM$^*$ model by Lemma \ref{L02}.
\medskip

Thus, $t_h^e(BH_n)\geq 2n-h$ under the MM$^*$ model for $0\leq h \leq 2n$ and $n\geq2$. $\qed$

Note that under the PMC model, Case 2 of Lemma \ref{LLL} is non-existent and the proof of Case 1 of Lemma \ref{LLL} also holds for $n\geq2$. Therefore, by Corollary \ref{LLU} and Lemma \ref{LLL}, Theorem \ref{main2} holds.
\begin{thm}\label{main2}
Let $BH_n$ be an $n$-dimensional balanced hypercube with $n\geq2$. Then $t_h^e(BH_n)=2n-h$ under the PMC model and the MM$^*$ model for $0\leq h \leq 2n$.
\end{thm}

Now, we determine the $h$-edge tolerable diagnosability of balanced hypercube $BH_1$.

\begin{thm}\label{r}
Let $BH_1$ be a $1$-dimensional balanced hypercube. Then
\begin{equation*}\label{}
t_h^e(BH_1)=
\left\{
  \begin{array}{lll}
   1, & \hbox{if $0\leq h \leq 1$}; \\
     0, & \hbox{if $h=2$}
  \end{array}
\right.
\end{equation*}
under the PMC model and $t_h^e(BH_1)=0$ under the MM$^*$ model for $0\leq h \leq 2$.
\end{thm}
\proof Note that $BH_1$ is isomorphic to a cycle with four vertices. Suppose $V(BH_1)=\{0, 1, 2, 3\}$ and $E(BH_1)=\{01, 12, 23, 30\}$ (see Figure \ref{DD1}). Since $BH_1$ is $2$-regular, $t_2^e(BH_1)=0$ under both diagnosis models.

Let $F_1=\{0, 1\}$, $F_2=\{2, 3\}$ and $F_e=\emptyset$. Then $F_1$ and $F_2$ are indistinguishable in $BH_1-F_e$ under the PMC model. Thus, $t_1^e(BH_1)\leq1$ and $t_0^e(BH_1)\leq 1$ under the PMC model.

On the other hand, for an arbitrary edge subset $F_e\subseteq E(BH_1)$ with $|F_e|\leq 1$, suppose that two distinct vertex subsets $F_1, F_2\subseteq V(BH_1)$ satisfy that $|F_1|\leq1$ and $|F_2|\leq1$. Then $F_1\cap F_2=\emptyset$ and $1\leq|F_1\bigtriangleup F_2|=|F_1\cup F_2|\leq2$. Since $BH_1-F_e$ is connected for $|F_e|\leq 1$, $N_{BH_1-F_e}(F_1\bigtriangleup F_2)=N_{BH_1-F_e}(F_1\cup F_2)\neq\emptyset$. Hence, there is an edge between $F_1\bigtriangleup F_2$ and $V(BH_1-F_e)-(F_1\cup F_2)$. By Lemma \ref{L01}, $t_h^e(BH_1)\geq 1$ for $0\leq h\leq1$ under the PMC model.

Let $F_1=\{0\}$, $F_2=\{2\}$ and $F_e=\emptyset$. Then $F_1$ and $F_2$ are indistinguishable in $BH_1-F_e$ under the MM$^*$ model. Thus, $t_h^e(BH_1)\leq 0$ for $0\leq h\leq1$ under the MM$^*$ model. Hence, $t_1^e(BH_1)=t_0^e(BH_1)=0$ under the MM$^*$ model.

This completes the proof of Theorem \ref{r}. $\qed$ 

\section{Conclusions}\label{4}
In this paper, we determine the $h$-edge tolerable diagnosability of balanced hypercubes $BH_n$ under the PMC model and the MM$^*$ model for $0\leq h \leq 2n$ and $n\geq1$ (see Table \ref{tab:2}). In particular, the traditional diagnosability of $BH_n$ is determined. Our future research interest is to investigate the $h$-edge tolerable diagnosability of a regular graph with triangles, which will provide a more precise measure for the fault diagnosis capability of a multiprocessor system.
\begin{table}[htbp]\normalsize
\begin{center}
\extrarowheight=4pt
\renewcommand{\arraystretch}{1.3}
\begin{tabular}{|c|c|c|c|}\hline
{\multirow{2}{*}{ }}& \multicolumn{2}{|c|}{~{\rm PMC} model~} &~{\rm MM}$^*${\rm ~model}~\\
\cline{2-4}
& $0\leq h\leq 2n-1$ & $h=2n$ & $0\leq h\leq 2n$\\
\hline
$n=1$ & $1$ & $0$ & $0$\\
\hline
$n\geq2$ & \multicolumn{2}{|c|}{$2n-h$} & $2n-h$\\
\hline
\end{tabular}
\end{center}
\caption{\quad \label{tab:2} The $h$-edge tolerable diagnosability of $BH_n$}

\end{table}

\noindent


\begin{thebibliography}{99}

\bibitem{Bon} J.A. Bondy, U.S.R. Murty, Graph Theory with Applications, The Macmillan Press Ltd, New York, 1976.

\bibitem{CDH} N.W. Chang, W.H. Deng, S.Y. Hsieh, Conditional diagnosability of (n, k)-star networks under the comparison diagnosis model, {\em IEEE Trans. Reliab.}, {\bf 64} (1) (2015), 132--143.

\bibitem{Dah} A.T. Dahbura, G.M. Masson, An $O(n^{2.5})$ faulty identification algorithm for diagnosable systems, {\em IEEE Trans. Comput.},  {\bf 33} (6) (1984), 486--492.

\bibitem{FanL} J. Fan, X. Lin, The $t/k$-diagnosability of the BC graphs, {\em IEEE Trans. Comput.},  {\bf 54} (2005), 176--184.

\bibitem{Gu} M.M. Gu, R.X. Hao, D.X. Yang, A short note on the $1, 2$-good-neighbor diagnosability of balanced hypercubes,  {\em J. Interconnect. Netw.}, {\bf 16} (2) (2016), pp. 1650001 (12 pages).

\bibitem{Lai} P.L. Lai, J.J.M. Tan, C.P. Chang, L.H. Hsu, Conditional diagnosability measures for large multiprocessor systems, {\em IEEE Trans. Comput.}, {\bf 54} (2) (2005), 165--175.

\bibitem{Lin17} L. Lin, L. Xu, R. Chen, S.Y. Hsieh, D.Wang, Relating extra connectivity and extra conditional diagnosability in
regular networks, {\em IEEE Trans. Depend. Secure Comput.} (2017), doi: 10.1109/TDSC.2017.2726541.

\bibitem{Mae} J. Maeng, M. Malek, A comparison connection assignment for self-diagnosis of multiprocessor systems, in: {\em  Proceeding of 11th International Symposium on Fault-Tolerant Computing}, 1981, pp. 173--175.

\bibitem{Pen} S.L. Peng, C.K. Lin, J.J.M. Tan, L.H. Hsu, The $g$-good-neighbor conditional diagnosability of hypercube under the PMC model, {\em Appl. Math. Comput.}, {\bf 218} (21) (2012), 10406--10412.

\bibitem{Pre} F.P. Preparata, G. Metze, R.T. Chien, On the connection assignment problem of diagnosis systems, {\em IEEE Trans. Electron. Comput.}, {\bf EC-16} (6) (1967), 848--854.

\bibitem{Sen} A. Sengupta, A. Dahbura, On self-diagnosable multiprocessor system: diagnosis by the comparison approach, {\em IEEE Trans. Comput.}, {\bf 41} (11) (1992), 1386--1396.

\bibitem{WanH} S. Wang, W. Han, The $g$-good-neighbor conditional diagnosability of $n$-dimensional hypercubes under the MM$^*$ Model, {\em Inform. Process. Lett.}, {\bf 116} (2016), 574--577.

\bibitem{Wei} Y. Wei, M. Xu, On $g$-good-neighbor conditional diagnosability of $(n, k)$-star networks, {\em Theoret. Comput. Sci.}, {\bf 697} (2017), 79--90.

\bibitem{Wei186} Y. Wei, M. Xu,  The $g$-good-neighbor conditional diagnosability of locally twisted cubes, {\em J. Oper. Res. Soc. China}, {\bf 6} (2) (2018), 333--347.

\bibitem{Wei18} Y. Wei, M. Xu, The $1, 2$-good-neighbor conditional diagnosabilities of regular graphs, {\em Appl. Math. Comput.}, {\bf 334} (2018), 295--310.

\bibitem{Wei188} Y. Wei, M. Xu, Hybrid fault diagnosis capability analysis of regular graphs, {\em Theoret. Comput. Sci.}, {\bf 760} (2019), 1--14.

\bibitem{WuH} J. Wu, K. Huang, The balanced hypercube: a cube-based system for fault-tolerant applications, {\em IEEE Trans. Comput.}, {\bf 46} (4) (1997), 484--490.

\bibitem{XuTH} M. Xu, K. Thulasiraman, X.D. Hu, Conditional diagnosability of matching composition networks under the PMC model, {\em IEEE Trans. Circuits Syst.}, II, Express Briefs {\bf 56} (11) (2009), 875--879.

\bibitem{XuTZ}  M. Xu, K. Thulasiraman, Q. Zhu, Conditional diagnosability of a class of matching composition networks under the comparison model, {\em Theoret. Comput. Sci.}, {\bf 674} (2017), 43--52.

\bibitem{Yang} M.C. Yang, Conditional diagnosability of balanced hypercubes under the PMC model, {\em Inform. Sci.}, {\bf 222} (2013), 754--760.

\bibitem{Yang13} M.C. Yang, Conditional diagnosability of balanced hypercubes under the MM$^*$ model, {\em J. Supercomput.}, {\bf 65} (3) (2013), 1264--1278.

\bibitem{YuanL} J. Yuan, A.X. Liu, X. Ma, X. Qin, J. Zhang, The $g$-good-neighbor conditional diagnosability of $k$-ary $n$-cubes under the PMC model and MM$^*$ model, {\em IEEE Trans. Parallel Distrib. Syst.}, {\bf 26} (4) (2015), 1165--1177.

\bibitem{Zha} S. Zhang, W. Yang, The $g$-extra conditional diagnosability and sequential $t/k$-diagnosability of hypercubes, {\em Int. J. Comput. Math.}, {\bf 93} (3) (2016), 482--497.

\bibitem{ZhaX} X. Zhang, L. Xu, L. Lin, Y. Huang, X. Wang, The $(t, k)$-diagnosability of balanced hypercube under the PMC model, {\em Int. J. Comput. Math. Comput. Syst. Theory}, {\bf 3} (4) (2018), 230--243.

\bibitem{Zho} S. Zhou, The conditional fault diagnosability of $(n, k)$-star graphs, {\em Appl. Math. Comput.}, {\bf 218} (19) (2012), 9742--9749.

\bibitem{ZhuL} Q. Zhu, L. Li, S. Liu, X. Zhang, Hybrid fault diagnosis capability analysis of
hypercubes under the PMC model and MM$^*$ model, {\em Theoret. Comput. Sci.}, {\bf 758} (2019), 1--8.
\end{thebibliography}
\end{document}